\documentclass[oneside]{amsart}
\usepackage{amsmath,graphicx}

\newtheorem{thm}{Theorem}[section]
\newtheorem{lem}[thm]{Lemma}

\theoremstyle{definition}
\newtheorem{defn}[thm]{Definition}

\newtheorem{conven}[thm]{Convention}

\newcommand{\blackboard}[1]{\ensuremath{\mathbb{#1}}}
\newcommand{\script}[1]{\ensuremath{\mathcal{#1}}}
\newcommand{\smallcaps}[1]{\textrm{\textsc{#1}}}
\newcommand{\mathbold}[1]{\ensuremath{\mathbf{#1}}}
\newcommand{\Z}{\blackboard{Z}}
\newcommand{\R}{\blackboard{R}}
\newcommand{\disc}{\mathbold{D}}
\newcommand{\braid}{\smallcaps{Braid}}
\newcommand{\pure}{\smallcaps{PBraid}}
\newcommand{\symm}{\smallcaps{Sym}}
\newcommand{\modd}{\smallcaps{Mod}}

\begin{document}

\title[Pure braids]{Geometric presentations for\\ the pure braid
  group}

\author[Margalit]{Dan Margalit}
  \address{Dept. of Mathematics\\
    University of Utah\\
    Salt Lake City, UT 84112-0090}
  \email{margalit@math.utah.edu}

\author[McCammond]{Jon McCammond}
  \address{Dept. of Mathematics\\
    U.C. Santa Barbara\\
    Santa Barbara, CA 93106-3080}
  \email{jon.mccammond@math.ucsb.edu}

\date{\today}

\begin{abstract}
We give several new positive finite presentations for the pure braid
group that are easy to remember and simple in form.  All of our
presentations involve a metric on the punctured disc so that the
punctures are arranged ``convexly'', which is why we describe them as
geometric presentations.  Motivated by a presentation for the full
braid group that we call the ``rotation presentation'', we introduce
presentations for the pure braid group that we call the ``twist
presentation'' and the ``swing presentation''.  From the point of view
of mapping class groups, the swing presentation can be interpreted as
stating that the pure braid group is generated by a finite number of
Dehn twists and that the only relations needed are the disjointness
relation and the lantern relation.
\end{abstract}

\maketitle

The braid group has had a standard presentation on a minimal
generating set ever since it was first defined by Emil Artin in the
1920s \cite{Ar26}.  In 1998, Birman, Ko, and Lee \cite{BiKoLe98} gave
a more symmetrical presentation for the braid group on a larger
generating set that has become fashionable of late (see, for example,
\cite{Be03}, \cite{Br01}, \cite{BrWa02}, or \cite{Kr00}).  Our goal is
to apply a similar idea to the pure braid group.  The standard finite
presentation for the pure braid group (also due to Artin \cite{Ar47})
is slightly complicated and not that easy to remember.  The
presentations introduced here are, we believe, simple, easy to
remember and intuitively clear.  The article is structured as follows.
In \S\ref{sec:braid} we present a variation of the Birman--Ko--Lee
presentation for the full braid group that we call the rotation
presentation, and in sections \ref{sec:pure}, \ref{sec:twist}, and
\ref{sec:swing} we establish increasingly simple presentations for the
pure braid group that we call the modified Artin presentation, the
twist presentation, and the swing presentation.  For these
presentations, we think of the braid group as the fundamental group of
the configuration space of $n$ points in the disk.  If we reinterpret
the swing presentation in terms of mapping class groups, we get a
presentation where the generators are Dehn twists, and the relations
are the disjointness relation and the lantern relation (see
section~\ref{sec:swing}).  The final section explores some possible
extensions.

\section{Braids}\label{sec:braid}

This section gives an unusual presentation of the full braid group
using the notion of a convexly punctured disc.  In addition to proving
that it is equivalent to the (closely related) Birman--Ko--Lee
presentation, we introduce several notions that pave the way for our
new presentations of the pure braid group.

The $n$-string braid group $\braid_n$ can be viewed as the fundamental
group of the configuration space of $n$ distinct but indistinguishable
points in a disc: $\braid_n \cong \pi_1(\mathcal{C}(\disc,n))$.  The points are
called \emph{punctures} and the elements of $\braid_n$ can be thought
of as homotopy classes of based loops in $\mathcal{C}(\disc,n)$, or equivalence
classes of \emph{motions} of these $n$ points in $\disc$ that start
and end at the same configuration.  The group gets its name from
Artin's original definition, which is different; the survey
\cite{BirBre04} proves these definitions are equivalent.  We begin by
defining some simple elements in $\braid_n$.

\begin{defn}[Half-twists]
Fix a configuration of the $n$ points in $\disc$ to serve as the
basepoint of $\pi_1(\mathcal{C}(\disc,n))$ and choose an embedded arc between
two punctures that only meets the set of punctures at its endpoints.
Let $\disc'$ be any subdisc of $\disc$ that contains the arc, the two
punctures it connects, and no other punctures.  Since
$\pi_1(\mathcal{C}(\disc',2)) \cong \mathbb{Z}$, we can define a
\emph{half-twist} to be a generator of $\pi_1(\mathcal{C}(\disc',2))$.  The
\emph{positive half-twist} corresponds to the half-twist where the two
points move around each other in a clockwise fashion.
\end{defn}

The braid group is generated by various finite sets of positive
half-twists.  In fact, any set of $n-1$ positive half-twists along
non-crossing arcs that connect the punctures in a tree-like fashion is
sufficient.  Typically people choose as basepoint the configuration
where the punctures lie in a straight line and they use the straight
arcs connecting neighboring punctures to define a generating set of
positive half-twists, but this is just a convenient standardization.

A more symmetric generating set is obtained by arranging the punctures
at the vertices of a convex $n$-gon in the disc, and using all of the
positive half-twists along the line segments connecting pairs of
punctures.  Notice that to make such a definition, the punctured disc
needs to be more than just a topological punctured disc: a metric
needs to be imposed so that convexity makes sense.

\begin{figure}[ht]
\begin{tabular}{ccc}
\includegraphics[scale=1]{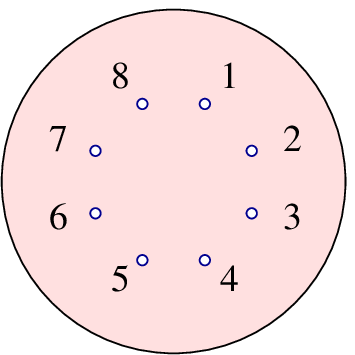} & \hspace*{2em} &
\includegraphics[scale=1]{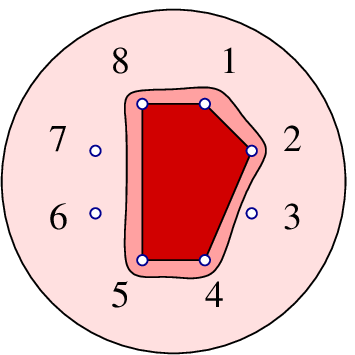}
\end{tabular}
\caption{A convexly punctured disc with $8$ punctures and a standard
labeling and the convexly punctured subdisc $\disc_B$ inside
$\disc_{[8]}$ when $B=\{1,2,4,5,8\}$.\label{fig:subdisc}}
\end{figure}

\begin{defn}[Convexly punctured discs]
Let $\disc$ be a topological disc in the Euclidean plane and assume
that $\disc$ has a distinguished $n$-element subset that we call its
\emph{punctures}.  If the disc is a convex subset of $\R^2$ and the
boundary of the convex hull of the set of punctures is an $n$-gon
(i.e. every puncture occurs as a vertex of the convex hull of the set
of punctures) then we say that $\disc$ is a \emph{convexly punctured
  disc} and that the punctures are in \emph{convex position}.  Let $P$
denote the convex $n$-gon whose vertices are the punctures.  There is
a natural cyclic ordering of the punctures corresponding to the
clockwise orientation of the boundary cycle of $P$.  A labeling of the
punctures is said to be \emph{standard} if it uses the set $[n] =
\{1,2,\ldots, n\}$ (or better yet $\Z/n\Z$) and the punctures are
labeled in the natural cyclic order.  See the left hand side of
Figure~\ref{fig:subdisc}.  More generally, when the punctures are
bijectively labeled by a finite set $A$, we refer to the convexly
punctured disc $\disc_A$.
\end{defn}

The notion of a convexly punctured disc is inherently recursive.  The
convexity of the arrangement of the punctures implies the existence of
a canonical convexly punctured subdisc corresponding to each nonempty
subset of punctures (see the right hand side of
Figure~\ref{fig:subdisc}).

\begin{defn}[Convexly punctured subdiscs]
Let $\disc_A$ be a convexly punctured disc.  If $B$ is any subset of
$A$, then there is a convexly punctured subdisc containing only the
punctures labeled by $B$.  In particular, define $\disc_B$ as an
$\epsilon$-neighborhood of the convex hull of the punctures labeled by
$B$.  Since the punctures in $\disc_A$ are in convex position, we can
choose $\epsilon$ small enough so that the only punctures in the
$\epsilon$-neighborhood are those labeled by $B$.  The resulting disc
is convex and the punctures are in convex position.
\end{defn}

\begin{figure}[ht]
\includegraphics[scale=1]{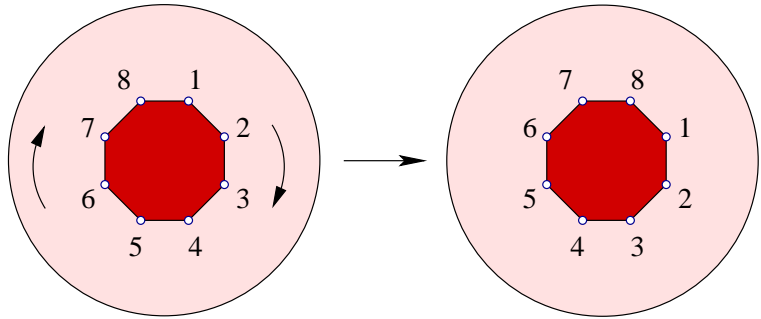}
\caption{Rotating all the punctures.\label{fig:rotate}}
\end{figure}

\begin{defn}[Rotating the punctures]
If $\disc_A$ is a convexly punctured disc, then there is a special
element of $\pi_1(\mathcal{C}(\disc_A,|A|))$ whose definition uses the convex
position of the punctures.  Let $P_A$ denote the convex hull of the
punctures and let $R_A$ be the motion which simultaneously moves each
puncture in $\disc_A$ clockwise along one side of the polygon $P_A$.
This is called \emph{rotating the punctures}.  See
Figure~\ref{fig:rotate}.  In the special case where $|A|=2$, the
motion we intend is the positive half-twist along the line segment
$P_A$, and when $|A|=1$ the rotation is the trivial motion.
\end{defn}

\begin{figure}[ht]
\includegraphics{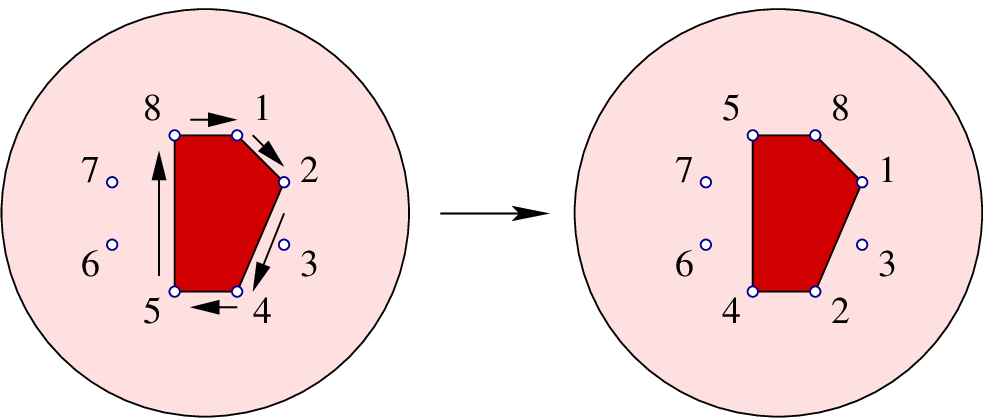}
\caption{The rotation $R_B$ inside $\disc_{[8]}$ when
  $B=\{1,2,4,5,8\}$.\label{fig:r-a}}
\end{figure}

Because of the recursive nature of convexly punctured discs, for each
$B \subset A$ there is a well-defined element $R_B$ inside
$\pi_1(\mathcal{C}(\disc_A,|A|))$ which rotates the punctures inside the subdisc
$\disc_B$ while leaving the remaining punctures fixed.  See
Figure~\ref{fig:r-a} for an illustration.

When we wish to emphasize the fact that the metric on $\disc_A$ is
used to define $R_B$, we call this a \emph{convex rotation}.  The
collection of all non-trivial convex rotations inside $\disc_A$,
denoted $\script{R}_A$, is the set $\{ R_B\ |\ B\subset A \textrm{
  with }|B| \geq 2\}$.  There are two types of relations among these
rotations that are easy to establish, but in order to state these
relations cleanly we need a pair of definitions.

\begin{figure}[th]
\begin{tabular}{ccc}
\includegraphics[scale=1]{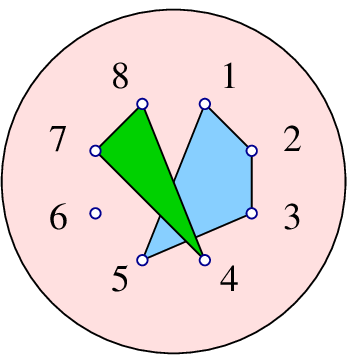} & \hspace*{2em} &
\includegraphics[scale=1]{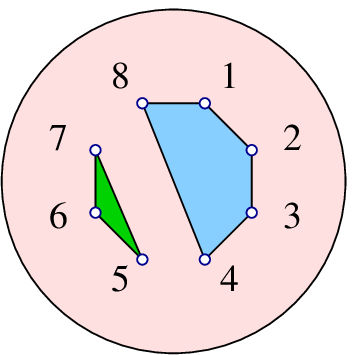}
\end{tabular}
\caption{The subsets $\{1,2,3,5\}$ and $\{4,7,8\}$ are crossing
subsets of $[8]$ and the subsets $\{1,2,3,4,8\}$ and $\{5,6,7\}$ are
non-crossing.\label{fig:crossing}}
\end{figure}

\begin{defn}[Non-crossing]
Let $\disc_A$ be a convexly punctured disc and let $B$ and $C$ be
disjoint subsets of $A$.  When the convex hull of $B$ and the convex
hull of $C$ do not intersect, $B$ and $C$ are said to be
\emph{non-crossing}. See Figure~\ref{fig:crossing}.  More generally,
an unordered collection $\{B_1,B_2,\ldots,B_k\}$ of pairwise disjoint
subsets of $A$ is called is non-crossing if $B_i$ and $B_j$ are
non-crossing for each $i\neq j$.
\end{defn}

\begin{figure}[th]
\begin{tabular}{ccc}
\includegraphics[scale=1]{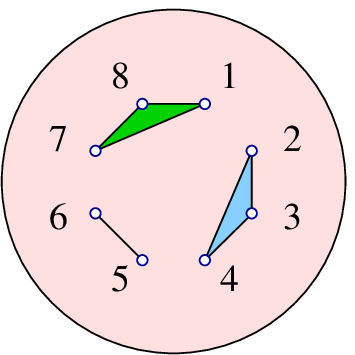} & \hspace*{2em} &
\includegraphics[scale=1]{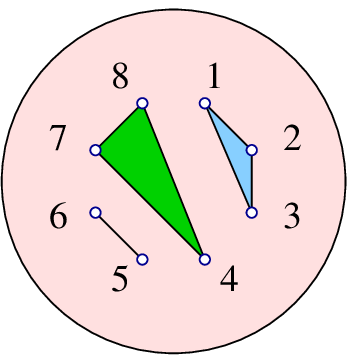}
\end{tabular}
\caption{If $B = \{2,3,4\}$, $C=\{5,6\}$ and $D=\{7,8,1\}$ then
  $(B,C,D)$ is an admissible partition, as is $(C,D,B)$, but $(C,B,D)$
  is not.  On the other hand, although the subsets $\{1,2,3\}$,
  $\{4,7,8\}$, and $\{5,6\}$ are non-crossing, no ordering of these
  three subsets is admissible.\label{fig:admissible}}
\end{figure}

\begin{defn}[Admissible partitions]
Let $\disc_A$ be a convexly punctured disc.  An ordered partition
$(A_1,A_2,\ldots,A_k)$ of $A$ is called an \emph{admissible partition}
(of $A$) if the cyclic ordering of the elements in $A$ is consistent
with the partial cyclic ordering determined by ordering of the $A_i$.
In other words, $(A_1,A_2,\ldots,A_k)$ is an admissible partition if
there is a point $x$ on the boundary of the convex hull of the
punctures so that the order in which the punctures occur starting at
$x$ and reading clockwise around the boundary consists of all of the
punctures labeled by $A_1$, followed by all of the punctures labeled
by $A_2$, and so on; see Figure~\ref{fig:admissible}. More generally,
$(B_1,B_2,\ldots, B_k)$ is \emph{admissible} inside $\disc_A$ if it is
an admissible ordering of the punctures inside the subdisc $\disc_B$
where $B= \cup_{i=1}^k B_i$.
\end{defn}

\begin{conven}
We follow the convention that uppercase letters such as $A$, $B$, and
$C$ denote sets while lowercase letters such as $i$, $j$, and $k$
denote elements.  For simplicity, unions of sets will be replaced by
juxtapositions and brackets around singleton sets will be removed when
there is no danger of confusion.  Thus $\{i\} \cup B \cup C$ will be
abbreviated as $iBC$.  Finally, we follow the standard practice in the
braid group literature \cite{BirBre04} and compose our motions from
left to right.  Thus, when we write $U\ V$, the motion denoted $U$
occurs first, followed by the motion denoted $V$.
\end{conven}

\begin{figure}[th]
\includegraphics[scale=1]{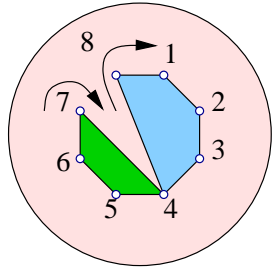}
\caption{An illustration that $R_{45678123}$ is the same as $ R_{4567}
\  R_{48123}$.\label{fig:r-split}}
\end{figure}

\begin{defn}[Rotation relations]\label{def:rotate-rel}
Let $\disc_A$ be convexly punctured disc, let $\braid_A$ be
$\pi_1(\mathcal{C}(\disc_A,|A|))$, and let $\script{R}_A \subset \braid_A$ denote
the set of convex rotations inside $\disc_A$.  The elements in
$\script{R}_A$ satisfy the following two types of relations.
\begin{eqnarray*}
R_B R_C = R_C R_B && \textrm{ when $B$ and $C$ are
  non-crossing}\\ R_{iBC} = R_{iB} R_{iC} && \textrm{ when
  $(\{i\},B,C)$ is admissible}
\end{eqnarray*}
The first type of relation holds because the rotations are occurring
in the disjoint subdiscs $\disc_B$ and $\disc_C$, and the second type
of relation is simply a factorization of the rotation into two smaller
rotations.  See Figure~\ref{fig:r-split} for an illustration.
\end{defn}

In Theorem~\ref{thm:braid} we show that the convex rotations subject
only to the rotation relations give a finite presentation for the
braid group.  To facilitate the proof we first review the
Birman--Ko--Lee presentation of $\braid_A$.

\begin{defn}[The Birman--Ko--Lee presentation]
Let $\disc_A$ be a convexly punctured disc and let $\braid_A$ be
$\pi_1(\mathcal{C}(\disc_A,|A|))$.  The presentation of $\braid_A$ introduced by
Birman, Ko, and Lee can be readily restated using the language we have
introduced above.  The generators they use are the rotations $R_{ij}$
and the presentation they give is the following:
\[ \braid_A \cong \left< \{ R_{ij}\}\  \begin{array}{|cl}
R_{ij}\  R_{kl} = R_{kl}\  R_{ij} & \textrm{ when $ij$ and
$kl$ are non-crossing}\\ R_{ijk} = R_{ij}\  R_{ik} & \textrm{ when
$(\{i\},\{j\},\{k\})$ is admissible} \end{array} \right>\] Strictly
speaking, they never introduce the elements $R_{ijk}$ and instead
equate its three possible factorizations, $R_{ij}\ R_{ik} =
R_{ik}\  R_{jk} = R_{jk}\  R_{ij}$, but the net effect is the
same.
\end{defn}

It is now easy to show that the convex rotations subject only to the
rotation relations given another presentation of the braid group.

\begin{thm}[Rotation presentation]\label{thm:braid}
If $\disc_A$ is a convexly punctured disc, $\braid_A$ is
$\pi_1(\mathcal{C}(\disc_A,|A|))$, and $\script{R}_A$ is the set of convex
rotations inside $\disc_A$, then
\[ \braid_A \cong \left< \script{R}_A\  \begin{array}{|cl}
R_B\  R_C = R_C\  R_B & \textrm{ when $B$ and $C$ are
non-crossing}\\ R_{iBC} = R_{iB}\  R_{iC} & \textrm{ when
$(\{i\},B,C)$ is admissible} \end{array} \right>\]
\end{thm}

\begin{proof}
Let $G_A$ be the group defined by this presentation, let $B_A$ be the
group defined by the Birman--Ko--Lee presentation and let $\braid_A$
be $\pi_1(\mathcal{C}(\disc_A,|A|))$.  Since all of the relations of $G_A$ hold
in $\pi_1(\mathcal{C}(\disc_A,|A|))$, there is a group homomorphism $f:G_A\to
\braid_A$.  Similarly, all of the relations of the Birman--Ko--Lee
presentation are included as relations in $G_A$, so there is a group
homomorphism $g:B_A\to G_A$.  Because the composition $f\circ g$ is
nothing other than the standard isomorphism between $B_A$ and
$\braid_A$, we know that $g$ must be one-to-one.  Finally, the
factorization rules show that the rotations of the form $R_{ij}$ are
enough to generate $G_A$.  Thus $g$ is onto.  Since $g$ and $f\circ g$
are both isomorphisms, so is $f$, and all three groups are isomorphic.
\end{proof}

Theorem~\ref{thm:braid} can be summarized as follows: if the basepoint
of $\pi_1(\mathcal{C}(\disc,n))$ is a convexly punctured disc, then this group
is generated by the convex rotations and all of its relations are
consequences of the semi-obvious rotation relations: disjoint
rotations commute and larger rotations can be factored into smaller
ones.

\section{Pure braids}\label{sec:pure}

In this section we give a new finite positive presentation for the
pure braid group that is similar in many ways to the rotation
presentation for the full braid group (Theorem~\ref{thm:braid}).  If
$\disc_A$ is a convexly punctured disc, then the \emph{pure braid
group} $\pure_A$ is the subgroup of $\braid_A$ where the punctures
must return to their original positions.  In other words, $\pure_A$ is
the fundamental group of the configuration space of $|A|$ distinct and
distinguishable points in $\disc$.  Algebraically, $\pure_A$ is the
kernel of the natural projection $\braid_A \to \symm_A$ that forgets
everything about the motion except the induced permutation of the
punctures (where $\symm_A$, of course, denotes the group of
permutations of $|A|$ elements).  Because the group $\pure_A$ is a
finite-index subgroup of $\braid_A$, every element of $\braid_A$ has a
power that lies in $\pure_A$.  This idea produces our first class of
pure braid elements.

\begin{defn}[Swinging the punctures]
Let $\disc_A$ be a convexly punctured disc and let $R_B$ be one of its
convex rotations.  The smallest power of $R_B$ that lies in $\pure_A$
is the motion $S_B=(R_B)^{|B|}$.  We refer to this as \emph{swinging
the punctures} labeled $B$, swinging being an apt term for a vigorous
rotation.  To emphasize the fact that the metric on $\disc_A$ is used
to define $S_B$, we call it a \emph{convex swing}.
\end{defn}

The collection of all convex swings inside $\disc_A$ is the finite set
$\script{S}_A = \{S_B\ |\ B \subset A \textrm{ with } |B|\geq 1\}$.
The convex swings around single punctures, though trivial, are
included in this set to facilitate extensions in later sections.  In
this section, the focus is on the convex swings of the form $S_{ij}$,
which are known to generate the pure braid group.  In fact, these were
the elements Artin used in his original presentation.

\begin{defn}[Artin's presentation of the pure braid group]
Let $\disc_A$ be a convexly punctured disc with a standard labeling
and let $\pure_A$ be its pure braid group.  Artin's presentation for
$\pure_A$ is generated by the elements of the form $S_{ij}$ subject to
the following five types of relations:
\begin{equation*}
S_{rs}^{-1} S_{ij} S_{rs} =
\left\{ \begin{array}{ll}
S_{ij} & \textrm{ if } r< s< i < j\\
S_{ij} & \textrm{ if } i< r< s< j\\
S_{rj} S_{ij} S_{rj}^{-1}  & \textrm{ if } r< i=s < j\\
(S_{ij} S_{sj}) S_{ij} (S_{ij} S_{sj})^{-1} & \textrm{ if }r= i
< s < j\\
(S_{rj} S_{sj} S_{rj}^{-1} S_{sj}^{-1}) S_{ij} (S_{rj} S_{sj}
S_{rj}^{-1} S_{sj}^{-1})^{-1} & \textrm{ if } r<i<s<j
\end{array} \right.\label{eq:artin}
\end{equation*}
\end{defn}

This presentation enabled Artin to establish a normal form for the
pure braids and to use this normal form to prove that $\pure_A$ is
poly-free (i.e. iteratively constructed using extensions by free
groups), but it seems clear from the final paragraph of the article
that he felt it has its limitations:

\begin{quote}
\textit{Although it has been proved that every braid can be deformed
  into a similar normal form the writer is convinced that any attempt
  to carry this out on a living person would only lead to violent
  protests and discrimination against mathematics.  He would therefore
  discourage such an experiment.}\\
\hspace*{3.2in}Emil Artin \cite{Ar47}
\end{quote}

We reformulate Artin's presentation is in terms of crossing,
non-crossing, and admissible partitions and eliminate the need for a
standard labeling.

\begin{thm}[Artin's presentation, modified] \label{thm:artin-mod}
Let $\disc_A$ be a convexly punctured disc and let $\pure_A$ be its
pure braid group.  The group $\pure_A$ is generated by the convex
swings $S_{ij}$ and every relation is can be derived from the
following three types of relations (assume all indices are distinct):
\begin{enumerate}
\item $[S_{ij},S_{rs}]=1$ when $\{i,j\}$ and $\{r,s\}$ are non-crossing,
\item $[S_{ij},S_{js}S_{rs}S_{js}^{-1}] = 1$ when $\{r,s\}$ and
  $\{i,j\}$ cross in cyclic order $r,i,s,j$,
\item $S_{sj}S_{rs}S_{rj} = S_{rs}S_{rj}S_{sj} = S_{rj}S_{sj}S_{rs}$
  when $(r,s,j)$ is admissible.
\end{enumerate}
\end{thm}

\begin{figure}[ht]
\includegraphics[scale=.75]{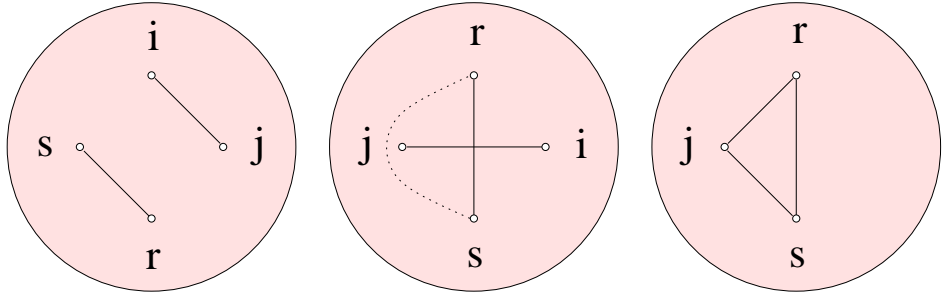}
\caption{The configurations of punctures for
  Theorem~\ref{thm:artin-mod}.  The pictures from left to right
  correspond to relations ($1$), ($2$), and
  ($3$).\label{fig:artin-mod}}
\end{figure}

Each of these relations can be viewed as an assertion that two
elements commute.  The configurations needed for each type of relation
are shown in Figure~\ref{fig:artin-mod}.  In relation ($2$), the
element $S_{js}S_{rs}S_{js}^{-1}$ corresponds to a (non-convex) swing
along the dotted arc.

\begin{proof}
Since the generators are the same and it is striaghtforward to check
that the given relations hold in $\pure_n$, it suffices to show that
Artin's original relations can be derived from relations ($1$), ($2$),
and ($3$). Relation ($1$) implies the first two of Artin's relations.
In Artin's third relation $i=s$.  After replacing $i$ with $s$ and
rearranging, the third relation is equivalent to $S_{sj}S_{rs}S_{rj} =
S_{rs}S_{rj}S_{sj}$ with $r<s<j$, which is the first equality of
relation ($3$).  In Artin's fourth relation $i=r$.  After replacing
$i$ with $r$ and rearranging, the fourth relation is equivalent to:
\begin{equation*}
\tag{$A4'$} S_{rj}S_{rs}S_{rj}S_{sj} = S_{rs}S_{rj}S_{sj}S_{rj}
\end{equation*}
\noindent
with $r<s<j$.  Relation ($A4'$) can be derived from the second equality
in relation ($3$) by starting with this relation, right multiplying both
sides by $S_{rj}$, and then applying the relation to the left hand
side as indicated:
\begin{eqnarray*}
S_{rj}S_{sj}S_{rs} &=& S_{rs}S_{rj}S_{sj} \\
S_{rj}(S_{sj}S_{rs}S_{rj}) &=& S_{rs}S_{rj}S_{sj}S_{rj} \\
S_{rj}S_{rs}S_{rj}S_{sj} &=& S_{rs}S_{rj}S_{sj}S_{rj}
\end{eqnarray*}
\noindent
Finally, Artin's fifth relation can be derived from relations ($2$)
and ($3$) as follows:
\begin{eqnarray*}
S_{ij}S_{sj}S_{rs}S_{sj}^{-1} &=&
S_{sj}S_{rs}S_{sj}^{-1}S_{ij} \\
S_{ij}(S_{sj}S_{rs}S_{rj})S_{rj}^{-1}S_{sj}^{-1} &=&
(S_{sj}S_{rs}S_{rj})S_{rj}^{-1}S_{sj}^{-1}S_{ij} \\
S_{ij}S_{rs}S_{rj}S_{sj}S_{rj}^{-1}S_{sj}^{-1} &=&
S_{rs}S_{rj}S_{sj}S_{rj}^{-1}S_{sj}^{-1}S_{ij}
\end{eqnarray*}
The first line is relation ($2$) slightly rearranged and the second
line freely reduces to the first.  We used relation ($3$) as indicated
to go from the second line to the third, which is Artin's fifth
relation rearranged.
\end{proof}

\section{Twists}\label{sec:twist}

Our second presentation of the pure braid group is defined in in terms
of ``convex twists'' that we like to think of as do-si-dos. A
\emph{do-si-do} is a movement in American square dancing where two
dancers approach each other and circle back to back (whence the French
term \emph{dos-\`a-dos}), and then return to their original positions.
The direction they face is unchanged throughout.  There is a similar
motion in the pure braid group.

\begin{figure}[ht]
\includegraphics[scale=.79]{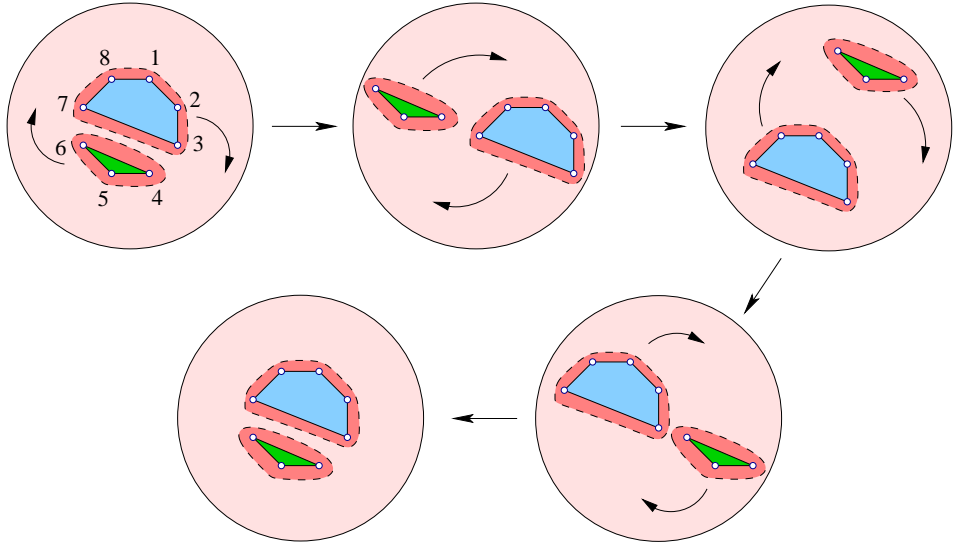}
\caption{The convex twist $T_{B,C}$ when $B=\{4,5,6\}$ and $C =
\{7,8,1,2,3\}$.\label{fig:dosido}}
\end{figure}

\begin{defn}[Twisting the punctures]
Let $\disc_A$ be a convexly punctured disc and let $B$ and $C$ be
non-crossing subsets that partition $A$.  Treating the convexly
punctured subdiscs $\disc_B$ and $\disc_C$ as though they were rigid
punctures, we can define a positive full twist between them.  In
keeping with `$R$' for rotation and `$S$' for swing, we use `$T$' for
twist and we denote this motion $T_{B,C}$.  The motion $T_{B,C}$ is
identical to $T_{C,B}$, so the subscripts should be considered
unordered.  As in the square dancing move, the subdiscs should be
moving by pure translations throughout. See Figure~\ref{fig:dosido}
for an illustration.  More generally, we define $T_{B,C}$ whenever $B$
and $C$ are merely non-crossing by having the twist take place inside
the convex subdisc $\disc_{BC}$. We call these elements \emph{convex
twists} and we let $\script{T}_A$ denote the collection of all convex
twists inside $\disc_A$ (i.e. $\script{T}_A = \{T_{B,C}\ |\ B,C
\subset A \textrm{ with $B$ and $C$ non-crossing}\}$).
\end{defn}

Notice that the element $T_{B,C}$ of $\pi_1(\mathcal{C}(\disc,n))$ can be
realized as a loop in $\mathcal{C}(\disc,n)$ where the points of $C$ stay fixed
throughout and the points of $B$ move around those of $C$.  This
perspective makes the third twist relation
(Definition~\ref{def:twist-rel}) easier to understand.  The next step
is to make some elementary observations about the relations satisfied
by the convex twists. An extra definition makes the relations easier
to state.

\begin{defn}[Nested pairs]
Let $\disc_A$ be a convexly punctured disc and let $(B,C)$ and $(D,E)$
be admissible subsets of $A$, not necessarily disjoint.  We say that
$(B,C)$ and $(D,E)$ are \emph{nested} if one of the following four
conditions hold: $BC \subset D$, $BC \subset E$, $DE \subset B$ or $DE
\subset C$.  For example, if $\disc_A$ is a convexly punctured disc
with a standard labeling by $A=[8]$ and $B=\{7,8,1,2,3\}$,
$C=\{4,5,6\}$, $D=\{7,1\}$ and $E= \{2,3\}$, then $(B,C)$ and $(D,E)$
are nested because $D \cup E$ is a subset of $B$.
\end{defn}

\begin{defn}[Twist relations]\label{def:twist-rel}
Let $\disc_A$ be convexly punctured disc, let $\pure_A$ be its pure
braid group, and let $\script{T}_A \subset \pure_A$ be the finite set
of convex twists inside $\disc_A$.  The elements in $\script{T}_A$
satisfy the following three types of relations that we call the
\emph{convex twist relations}:
\begin{eqnarray*}
T_{B,C}\  T_{D,E} = T_{D,E}\  T_{B,C} && \textrm{ when $BC$ and
  $DE$ are non-crossing} \\
T_{B,C}\  T_{D,E} = T_{D,E}\  T_{B,C} && \textrm{ when $(B,C)$
  and $(D,E)$ are nested} \\
T_{B,CD} = T_{B,C}\  T_{B,D} && \textrm{ when $(B,C,D)$ is admissible}
\end{eqnarray*}
\end{defn}

\begin{figure}[ht]
\includegraphics[scale=1]{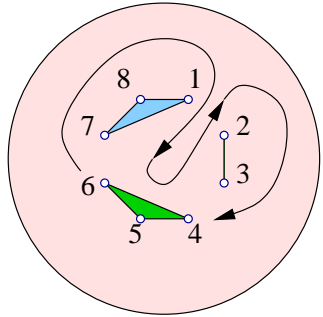}
\caption{An illustration that $T_{B,CD}$ is the same as $T_{B,C}\
T_{B,D}$ when $B=\{4,5,6\}$, $C=\{7,8,1\}$ and $D=\{2,3\}$.
\label{fig:t-split}}
\end{figure}

The first relation holds because the twists are occurring in the
disjoint subdiscs $\disc_{BC}$ and $\disc_{DE}$, and the third
relation is simply a factorization of the twist into two smaller
twists.  See Figure~\ref{fig:t-split} for an illustration.  Thus, the
only relation that need to be explained is the second one.  In this
case, it is useful to prove a stronger result first.

\begin{lem}[Twists and braids]\label{lem:twist-braid}
If $\disc_A$ is a convexly punctured disc and $T_{B,C}$ is a convex
twist in $\disc_A$, then $T_{B,C}\ U = U\ T_{B,C}$ for all $U$ in
$\braid_B$.
\end{lem}

\begin{figure}[ht]
\includegraphics[scale=1]{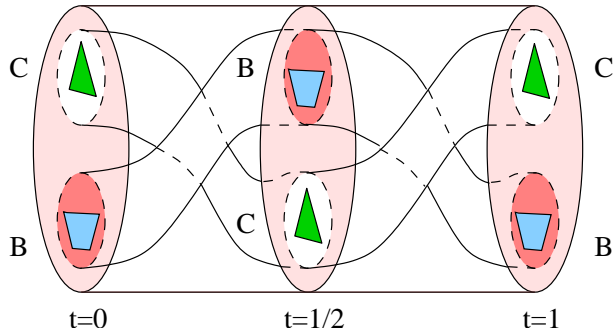}
\caption{The $3$-dimensional trace of a convex twist over
time.\label{fig:twist-braid}}
\end{figure}

\begin{proof}
The $3$-dimensional model of a braid (in Artin's original definition)
is obtained from its representation as an element of
$\pi_1(\mathcal{C}(\disc_A,|A|))$ by tracing the paths of the
punctures in the disc over time.  This gives a solid cylinder $\disc_A
\times [0,1]$ with $|A|$ strands inside it.  If we keep track of the
convex subdisc $\disc_B$ during the convex twist $T_{B,C}$, the result
looks something like Figure~\ref{fig:twist-braid}.  The key
observation is the solid tube which tracks $\disc_B$ over time is
internally untwisted.  Thus the action of any element of
$\pi_1(\mathcal{C}(\disc_B,|B|))$ on $\disc_B$ that takes place after
$T_{B,C}$ can be pushed back through this tube so that takes place
before $T_{B,C}$.
\end{proof}

The reason why $T_{B,C}$ and $T_{D,E}$ commute when $(B,C)$ and
$(D,E)$ are nested should now be clear.  If, for example, $DE \subset
B$, then $T_{D,E}$ is an element of $\braid_B$ and by
Lemma~\ref{lem:twist-braid} they commute.  The other three cases are
similar.  We show in Theorem~\ref{thm:twist} that the convex twists
subject only to the convex twist relations give a (finite positive)
presentation of the pure braid group.  To facilitate the proof we
first establish that the modified Artin relations
(Theorem~\ref{thm:artin-mod}) can be derived from the convex twist
relations, starting with the following lemma.

\begin{lem}\label{lem:triangles}
Let $\disc_A$ be a convexly punctured disc.  If $(B,C,D)$ is
admissible, then \[T_{C,B}\ T_{B,D}\ T_{D,C} = T_{B,D}
\ T_{D,C}\ T_{C,B} = T_{D,C}\ T_{C,B}\ T_{B,D}\] holds in $\pure_A$
and these relations are consequences of the convex twist relations.
\end{lem}

\begin{proof}
To prove that the first two expressions are equal it suffices to show
that $T_{C,B}$ commutes with $T_{B,D}\ T_{D,C}$. But $T_{B,D}\ T_{D,C}
= T_{D,B}\ T_{D,C} = T_{D,BC}$ by the factoring relation (since
$(D,B,C)$ is also admissible) and this commutes with $T_{C,B}$ because
nested twists commute.  The second equality is proved similarly.
\end{proof}

\begin{lem}\label{lem:twist-rel}
Let $\disc_A$ be a convexly punctured disc.  The modified Artin
relations (Theorem~\ref{thm:artin-mod}) are derivable from the convex
twist relations (Definition~\ref{def:twist-rel}).
\end{lem}

\begin{proof}
Since the convex swing $S_{ij}$ is another name for the convex twist
$T_{i,j}$, each modified Artin relation can be easily rewritten in
terms of the convex twists.  The first relation of
Theorem~\ref{thm:artin-mod} is covered by the assertion that
non-crossing convex twists commute and the third relation is a special
case of Lemma~\ref{lem:triangles}, so only the second relation remains
to be derived.  For later use, note that the second equality of
relation ($3$) in Theorem~\ref{thm:artin-mod} can be written as:
\begin{equation}\tag{$3'$} T_{s,j}\  T_{r,s}\  T_{s,j}^{-1} =
  T_{r,j}^{-1}\  T_{r,s}\  T_{r,j} 
\end{equation}
The fact that this relation holds in $\pure_A$ is clear from
Figure~\ref{fig:artin-mod}; both sides describe the nonconvex twist
along the dotted arc.  To derive the second modified Artin relation,
we start with a pair of commuting nested convex twists (second convex
twist relation), decompose them into smaller convex twists (third
convex twist relation), rearrange the equality, and finally apply
relation ($3'$) to the right hand side as follows:
\begin{eqnarray*}
T_{ris,j}\  T_{r,s} &=& T_{r,s}\  T_{ris,j} \\
T_{r,j}\  T_{i,j}\   T_{s,j}\  T_{r,s} &=& T_{r,s}\
T_{r,j}\  T_{i,j}\  T_{s,j}  \\
T_{i,j}\  T_{s,j}\  T_{r,s}\  T_{s,j}^{-1}  &=&
(T_{r,j}^{-1}\  T_{r,s}\  T_{r,j})\   T_{i,j} \\
T_{i,j}\   (T_{s,j}\  T_{r,s}\  T_{s,j}^{-1}) &=&
(T_{s,j}\  T_{r,s}\  T_{s,j}^{-1})\   T_{i,j} \\
\end{eqnarray*}
The last line is the second modified Artin relation.
\end{proof}

\begin{thm}[Twist presentation]\label{thm:twist}
If $\disc_A$ is a convexly punctured disc, then its pure braid group
is generated by convex twists and all of its relations are
consequences of the convex twist relations.  In particular, $\pure_A$
is isomorphic to the group defined by the following finite
presentation:
\[\left<\ \script{T}_A\ \begin{array}{|cl}
T_{B,C}\ T_{D,E} = T_{D,E}\ T_{B,C} & \textrm{ when $BC$ and $DE$ are
  non-crossing} \\ T_{B,C}\ T_{D,E} = T_{D,E}\ T_{B,C} & \textrm{ when
  $(B,C)$ and $(D,E)$ are nested} \\ T_{B,CD} = T_{B,C}\ T_{B,D} &
\textrm{ when $(B,C,D)$ is admissible}
\end{array} \right>
\]
\end{thm}

\begin{proof}
Let $G_A$ be the group defined by this listed presentation, let $PB_A$
be the group defined by Artin's presentation and let $\pure_A$ be the
pure braid group of $\disc_A$.  Since all of the relations of $G_A$
hold in the pure braid group $\pure_A$, there is a group homomorphism
$f:G_A\to \pure_A$.  Similarly, by Lemma~\ref{lem:twist-rel}, all of
the relations of Artin's presentation are induced by relations in
$G_A$, so there is a group homomorphism $g:PB_A\to G_A$.  Because the
composition $f\circ g$ is nothing other than the standard isomorphism
between $PB_A$ and $\pure_A$, we know that $g$ must be one-to-one.
Finally, the factorization rules show that the twists of the form
$T_{i,j}$ are enough to generate $G_A$.  Thus $g$ is also onto.  Since
$g$ and $f\circ g$ are both isomorphisms so is $f$, and all three
groups are isomorphic.
\end{proof}

Said differently, and perhaps more memorably, the pure braid group is
generated by do-si-dos and the only relations needed are: nested
do-si-dos commute, non-crossing do-si-dos commute, and do-si-dos
decompose.

\section{Swings}\label{sec:swing}

In this section we introduce our final finite presentation for the
pure braid group, this time inspired by mapping class groups and
particulary simple in form.  The \emph{relative mapping class group} of the
pair $(\Sigma,\Delta)$, where $\Sigma$ is a surface and $\Delta
\subset \Sigma$, is denoted $\smallcaps{Mod}(\Sigma,\Delta)$, and is
defined by
\[\smallcaps{Mod}(\Sigma,\Delta) = \pi_0(\textrm{Homeo}^+(\Sigma, \Delta))\]
where $\textrm{Homeo}^+(\Sigma,\Delta)$ denotes the set of orientation preserving
homeomorphisms of $\Sigma$ that fix the set $\Delta$ pointwise.  The
set $\modd(\Sigma,\Delta)$ forms a group under
function composition.  In order to be consistent with the earlier
sections, we continue to use the algebraic rather than functional
convention for these compositions (i.e., left to right, not right to
left).

Mapping class groups are relevant because the pure braid group can be
described in this language:
\[\pure_A \cong \smallcaps{Mod}(\disc_A,A \cup \partial \disc_A) \]
where $\disc_A$. For a proof of this isomorphism, see \cite{BirBre04}.

We begin by defining the simple elements of the mapping class
group that are traditionally used as generators.

\begin{figure}[ht]
\includegraphics{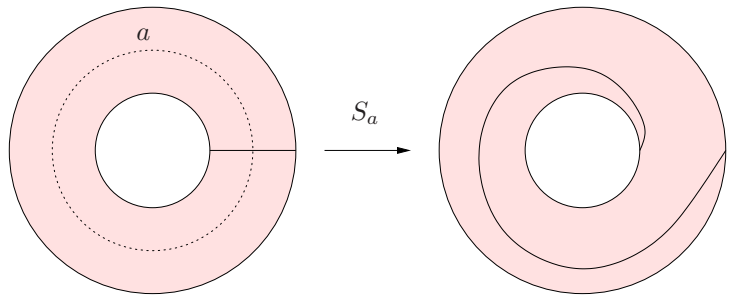}
\caption{The non-identity portion of a Dehn twist.\label{fig:dehn-twist}}
\end{figure}

\begin{defn}[Dehn twists]\label{def:dehn-twists}
If $\alpha$ is a simple closed curve in $\disc_A$ disjoint from $A
\cup \partial \disc_A$ then there is a homeomorphism of $\disc_A$ that
looks like Figure~\ref{fig:dehn-twist} on a small regular neighborhood
of $\alpha$ and the identity outside of this neighborhood.  The
homeomorphism described is well-defined up to an isotopy of $\alpha$,
so long as the intermediate curves remain disjoint from $A\cup
\partial \disc_A$.  If $a$ denotes this isotopy class of curves, the
element of $\pure_A$ it defines is called the \emph{Dehn twist} about
$a$ and denoted $S_a$.  The reason for this notation is discussed
below.
\end{defn}

\begin{defn}[Convex Dehn twists]\label{def:convex-dehn}
A simple closed curve $\alpha$ in a convexly punctured disc $\disc_A$
is \emph{convex} if its interior (i.e., the component of
$\disc-\alpha$ disjoint from $\partial \disc$) is convex.  We call the
Dehn twists about isotopy classes of convex simple closed curves
\emph{convex Dehn twists}.  Let $\mathcal{D}_A$ denote the set of all
convex Dehn twists in $\pure_A$.  This set is finite since each
isotopy class is uniquely determined by the subset $B$ of punctures
contained in the interior of any representative convex curve.  The
trivial class that surrounds no punctures is excluded from this set.
\end{defn}

Convex Dehn twists correspond, in fact, to motions we have already
encountered.

\begin{lem}[Convex swings and convex Dehn twists]\label{lem:dehn=swing}
Let $\disc_A$ be a convexly punctured disc.  Every convex swing $S_B$,
when viewed as an element of the mapping class group, is equal to the
convex Dehn twist $S_b$ where $b$ is the isotopy class of the boundary
cycle of the convex subdisc $\disc_B$.  Conversely, every convex Dehn
twist $S_b$ corresponds to the convex swing $S_B$ where $B$ labels the
set of punctures contained in the interior of $b$.
\end{lem}

It is because of this close connection between convex swings and
convex Dehn twists that we have chosen to use $S_b$ to denote a Dehn
twist even when the isotopy class $b$ does not contain a convex
representative.  In the sequel, we use the Dehn twist notation ($S_b$)
when no assumptions are made about convexity, but we typically switch
to swing notation ($S_B$) when all the isotopy classes contain convex
curves.  Dehn twists satisfy three well-known types of relations.

\begin{figure}[ht]
\includegraphics{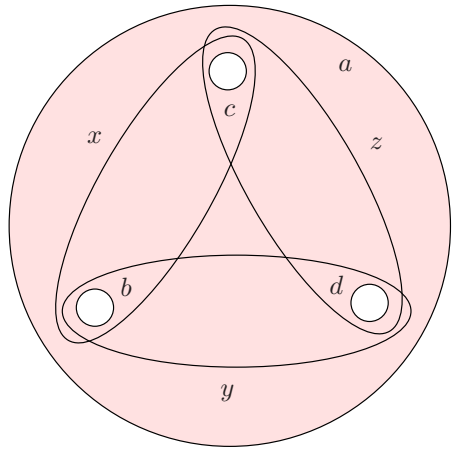}
\caption{The curves in the lantern relation.\label{fig:lantern}}
\end{figure}

\begin{defn}[Dehn twist relations]\label{def:dehn-twist-rel}
Let $\disc_A$ be punctured disc and let $\pure_A$ be its pure braid
group.  The Dehn twists in $\pure_A$ satisfy the following:
\[\begin{array}{rcll}
S_b &=& 1 & \textrm{ when $b$ surrounds a single puncture}\\
S_bS_c &=& S_cS_b & \textrm{ when $b$ and $c$ have disjoint
  representatives}\\
S_aS_bS_cS_d &=& S_xS_yS_z & \textrm{ when the representative curves
  look like Figure~\ref{fig:lantern}}
\end{array}\]
The first relation is trivial.  The second relation, called the
\emph{disjointness relation}, is obvious since the annuli in which the
twisting takes place can be chosen to be disjoint.  The third
relation, called the \emph{lantern relation}, was known to Dehn
\cite{md}.  The lantern relation is the only relation where the order
matters and only on one side.  If we were using functional notation,
the lantern relation would be $S_dS_cS_bS_a=S_zS_yS_x$.  The left hand
side could be rewritten as $S_aS_bS_cS_d$ since these Dehn twists
pairwise commute, but the order of $S_x$, $S_y$, and $S_z$ is
important.
\end{defn}

The convex twist version of the lantern relation is established in the
course of proving Lemma~\ref{lem:swing-rel}.  The first step is to
convert the Dehn twist relations to convex form.  Let $\disc_A$ be a
convexly punctured disc.  We call two sets $B,C\subset A$
\emph{compatible} if $B \subset C$, $C \subset B$, or $B$ and $C$ are
non-crossing.  This definition is useful since these are exactly the
conditions under which the boundary curves $\partial \disc_B$ and
$\partial \disc_C$ can be chosen to be disjoint.

\begin{defn}[Swing relations]\label{def:swing-rel}
If $\disc_A$ is a convexly punctured disc and we restrict our
attention to the Dehn twist relations that hold among the convex Dehn
twists, then we can rewrite them as relations among the convex swings.
We call these the \emph{convex swing relations}.
\[\begin{array}{rcll}
S_B &=& 1 & \textrm{ when $|B|=1$}\\ S_BS_C &=& S_CS_B & \textrm{ when
$B$ and $C$ are compatible} \\ S_{BCD}S_BS_CS_D &=&S_{CB}S_{BD}S_{DC}
& \textrm{ when $(B,C,D)$ is admissible}\\
\end{array}\]
Notice that if we choose isotopy classes of convex curves $a$, $b$,
$c$, $d$, $x$, $y$ and $z$ so that they surround the puncture sets
$BCD$, $B$, $C$, $D$, $BC$, $BD$, and $CD$, respectively, then these
seven curves can be arranged as in Figure~\ref{fig:lantern}.
\end{defn}

In Lemma~\ref{lem:swing-rel} we show that the convex twist relations
can be derived from the convex swing relations.  To show this we first
need to establish the connection between convex swings and convex
twists.

\begin{lem}[Twists as Swings]\label{lem:twist-to-swing}
Let $\disc_A$ be a convexly punctured disc.  Every convex twist
$T_{B,C}$ in $\disc_A$ can be rewritten as a product of three
commuting convex swings (or their inverses).  In particular, $T_{B,C}
= S_B^{-1}S_C^{-1} S_{BC}$ written in any order and $S_{BC} = S_B S_C
T_{B,C}$ written in any order.
\end{lem}

\begin{figure}[ht]
\includegraphics{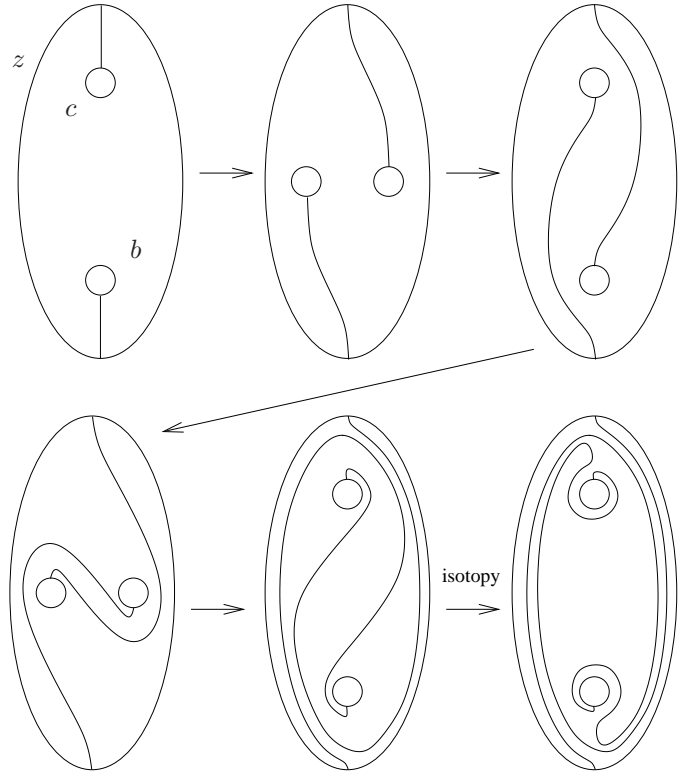}
\caption{A visual proof of
  Lemma~\ref{lem:twist-to-swing}.\label{fig:slow-twist}}
\end{figure}

\begin{proof}
One approach is to first establish the equality $S_{BC} = S_B S_C
T_{B,C}$.  The natural motion for $S_{BC}$ can be viewed as enacting
the motions $S_B$, $S_C$, and $T_{B,C}$ simultaneously with the convex
swings $S_B$ and $S_C$ taking place inside the convex subdiscs
$\disc_B$ and $\disc_C$ as they perform their do-si-do.  As in the
proof of Lemma~\ref{lem:twist-braid}, the motions $S_B$ and $S_C$ can
be pushed forward to back through these untwisted tubes so that the
resulting elements take place in any order.  Alternatively, let $b$,
$c$, and $z$ be the isotopy classes of the convex boundary cycles of
$\disc_B$, $\disc_C$ and $\disc_{BC}$, respectively.  A visual proof
that $T_{B,C}$ is equal to the product $S_b^{-1} S_c^{-1} S_z$ of
three commuting Dehn twists is shown in Figure~\ref{fig:slow-twist}.
\end{proof}

An easy consequence of Lemma~\ref{lem:twist-to-swing} is that the
swing $S_A$ is central in $\pure_A$.

\begin{thm}[$S_A$ is central]\label{thm:central}
If $\disc_A$ is a convexly punctured disc then for all $i,j\in A$,
there is an element $U$ such that $S_{ij} U = U S_{ij} = S_A$.  As a
consequence, the element $S_A$ is central in $\pure_A$.
\end{thm}

\begin{proof}
When $\{i,j\} = A$ we can set $U$ equal to the identity element and
there is nothing to prove, so assume $|A|>2$.  If $i$ and $j$ are
consecutive in the standard cyclic order (say $i$ followed by $j$)
then there is a nonempty set $B$ where $(i,j,B)$ is an admissible
partition of $A$.  By Lemma~\ref{lem:twist-to-swing}, $S_A = S_{ij}
(S_B T_{ij,B}) = (S_B T_{ij,B}) S_{ij}$, so the assertion is true with
$U = S_B T_{ij,B}$.  Finally, when $i$ and $j$ are not consecutive,
there are nonempty sets $B$ and $C$ such that $(i,B,j,C)$ is an
admissible partition of $A$.  Applying Lemma~\ref{lem:twist-to-swing}
twice we find that $S_A = S_{ijC} S_B T_{B,ijC} = (S_{ij} S_C
T_{ij,C}) S_B T_{B,ijC}$.  Since all five of these elements pairwise
commute, the assertion is true with $U = S_B S_C T_{ij,C} T_{B,ijC}$.
The final assertion is nearly immediate. If we select the element $U$
that is paired with the element $S_{ij}$, then parenthesizing $S_{ij}
U S_{ij}$ in two different ways shows that $S_A$ commutes with
$S_{ij}$ for all $i,j\in A$. Because these elements generate
$\pure_A$, $S_A$ is central.
\end{proof}

Theorem~\ref{thm:central} is, of course, obvious from the point of
view of mapping class groups (since disjoint twists commute),
configuration spaces (since nested swings commute), and in Artin's
original pictorial definition (just by staring at the pictures).  The
point here is that the formal presentation we are introducing makes it
easy to establish facts such as this without appealing to topological
intuition.

\begin{lem}[Swings as Twists]\label{lem:swing-to-twist}
Let $\disc_A$ be a convexly punctured disc.  Every convex swing $S_B$
in $\disc_A$ can be rewritten as a product of convex twists.
\end{lem}

\begin{proof}
Using the identity $S_{BC} = S_B S_C T_{B,C}$ from
Lemma~\ref{lem:twist-to-swing} repeatedly we can decompose any convex
swing involving more than $2$ punctures into a convex twist and two
convex swings involving strictly fewer punctures.  After finitely many
such steps the original convex swing has been rewritten as a product
of convex twists and convex swings involving single punctures.  Since
the latter are trivial, they drop out of the product, proving that the
original convex swing is a product of convex twists.
\end{proof}

\begin{lem}\label{lem:swing-rel}
Let $\disc_A$ be a convexly punctured disc.  The convex twist
relations (Theorem~\ref{thm:twist}) are derivable from the convex
swing relations (Definition~\ref{def:dehn-twist-rel}).
\end{lem}

\begin{proof}
Let $T_{B,C}$ and $T_{D,E}$ be nested or non-crossing convex twists.
When interpreted as elements of the mapping class group using
Lemma~\ref{lem:twist-to-swing}, both elements become a product of
three convex swings or their inverses.  The nested / non-crossing
hypothesis implies that all six convex curves can be chosen to be
simultaneously pairwise disjoint.  The fact that these convex twists
commute now follows immediately from the commuting convex swing
relations.

Finally, the third twist relation is closely related to the lantern
relation.  Let $(B,C,D)$ be an admissible partition of $\disc_A$.
Using Lemma~\ref{lem:twist-to-swing}, the third twist relation
$T_{B,CD} = T_{B,C} T_{B,D}$ is equivalent to $S_{BCD} S_B^{-1}
S_{CD}^{-1} = (S_{BC} S_B^{-1} S_C^{-1}) (S_{BD} S_B^{-1} S_D^{-1})$.
This is the relation we wish to establish.  Since the convex swings
$S_B$, $S_C$, and $S_D$ commute with all of the other convex swings in
the relation, they can be collected on the left hand side, cancelling
out an $S_B^{-1}$ in the process.  Finally, multiplying both sides on
the right by $S_{CD}$ produces the convex swing version of the lantern
relation.  Since these steps are reversible, the third twist relation
can be derived from the lantern relation and the disjointness
relation.
\end{proof}

\begin{thm}[Swing presentation]\label{thm:swing}
If $\disc_A$ is a convexly punctured disc, then its relative mapping
class group (i.e., its pure braid group) is generated by its convex
Dehn twists and all of its relations are consequences of the obvious
triviality and disjointness relations among the generators, together
with a finite number of lantern relations derived from admissible
partitions.  More concretely, $\pure_A$ is isomorphic to the group
defined by the following finite convex swing presentation:
\[\left<\ \script{S}_A\ \begin{array}{|cl}
S_B = 1 & \textrm{ when $|B|=1$}\\
S_B S_C = S_C S_B & \textrm{ when $B$ and $C$ are compatible} \\
S_{BCD} S_B S_C S_D =S_{CB} S_{BD} S_{DC} & \textrm{ when $(B,C,D)$
  is admissible} \\
\end{array} \right>\]
Alternatively (and equivalently), $\pure_A$ is isomorphic to the group
defined by the following finite convex Dehn twist presention:
\[\left<\ \script{S}_A\ \begin{array}{|cl}
S_b = 1 & \textrm{ when $b$ surrounds a single puncture}\\
S_b S_c = S_c S_b & \textrm{ when $b$ and $c$ have disjoint
representatives} \\ S_aS_bS_cS_d = S_zS_yS_x & \textrm{ when these
isotopy classes look like Figure~\ref{fig:lantern}}\\
\end{array} \right>\]
\end{thm}

\begin{proof}
Let $H_A$ be the group defined by the convex swing presentation and
let $G_A$ be the group defined by the convex twist presentation given
in Theorem~\ref{thm:twist}.  There is a group homomorphism $f:H_A \to
\pure_A$ since each of the relations in the presentation are known to
hold in the pure braid group, and there is a group homorphism $g:G_A
\to H_A$ extending the natural map that rewrites convex twists as a
product of convex swings since every convex twist relation in the
presentation of $G_A$ can be derived from the convex swing relations
in the presentation of $H_A$ (Lemma~\ref{lem:swing-rel}).  Because the
composition $f \circ g$ is the previously established isomorphism
between $G_A$ and $\pure_A$ (Theorem~\ref{thm:twist}), the map $g$ is
injective.  The map $g$ is also onto since by
Lemma~\ref{lem:swing-to-twist} the set $g(G_A)$ includes the
generating set of $H_A$.  Because $g$ and $f \circ g$ are
isomorphisms, so is $f$.  Finally, the conversion from the convex
swing presentation to the convex Dehn twist presentation is merely a
change of notation.
\end{proof}

\section{Final Comments}\label{sec:open}

In this final section we mention some possible extensions.  The (pure)
braid groups are paradigmatic examples of two distinct classes of
groups: mapping class groups and Artin groups.

\subsection{Other mapping class groups}

In the world of mapping class groups there is a natural modification
where the punctures are replaced by boundary components.

\begin{thm}
Let $\disc_A$ be a convexly punctured disc and let $\disc'_A$ be
$\disc_A$ with small open neighborhoods of the punctures removed.  The
relative mapping class group of $\disc'_A$ is generated by its convex
Dehn twists and all of its relations are consequences of the obvious
disjointness relations among the generators, together with a finite
number of lantern relations derived from admissible partitions.  More
concretely, $\modd(\disc'_A,\partial \disc'_A)$ is isomorphic to the
group defined by the following finite presentation:
\[\left<\ \script{S}_A\ \begin{array}{|cl}
S_B S_C = S_C S_B & \textrm{ when $B$ and $C$ are compatible} \\
S_B S_C S_D S_{BCD}=S_{CB} S_{BD} S_{DC} & \textrm{ when $(B,C,D)$
  is admissible} \\
\end{array} \right>\]
\end{thm}

The group defined above is an abelian extension the pure braid group
since there is a natural map from it to $\pure_A$ and the kernel is
generated by the elements $S_i$ which are central.  If $\disc''_A$ is
$\disc_A$ with small neighborhoods of only some of the punctures
removed then its relative mapping class group has a similar positiive
finite presentation with all of the convex disjointness and lantern
relations and only those triviality relations that correspond to the
remaining punctures.

The presentation for the relative mapping class group of a sphere with
discs removed naturally leads to nice presentations for the pure
stabilizers of any set of closed curves in the mapping class group of
a closed surface whose complements have genus zero and there is at
least a chance that a relatively simple presentation for the full
mapping class group of a closed surface could result.

\subsection{Other Artin groups}

In \cite{Al02} Daniel Allcock gave orbifold descriptions for most of
the other irreducible Artin groups of finite-type (i.e. those that
correspond to the irreducible finite Coxeter groups).  Because of the
existence of an orbifold description, there should be presentations
for the pure Artin groups of finite-type that are very similar in
nature to the ones given here.  In particular, for each pure Artin
group $G$ of finite-type the squares of the standard dual generators
should be a generating set, these generators should be identifiable
with some set of basic (convex) moves in the orbifold picture, and
there should be enough relations among the geometrically obvious
convex commutations, factorizations, and lantern relations to define a
presentation for $G$.

\vspace{2mm}
\noindent
\textit{Acknowledgements:} The authors would like to thank Angela
Barnhill and Saul Schleimer for helpful conversations and the National
Science Foundation for their financial support.  The first author was
partly supported by an NSF postdoctoral research fellowship and the
second author by an NSF standard grant.

\bibliography{refs-pure}
\bibliographystyle{plain}

\end{document}